\documentclass[12pt,oneside]{amsart}

\usepackage{amssymb}
\usepackage{rotating}

\usepackage{amsxtra}

\usepackage{pstricks}
\usepackage[all]{xy}
\usepackage{showidx}
\usepackage{amscd}
\usepackage[active]{srcltx}
\usepackage{multicol}
\usepackage{fancyhdr}
\usepackage{changebar}

\usepackage{amsfonts}
\usepackage{dsfont}
\usepackage{amsxtra}

\usepackage{index}

\theoremstyle{plain}

\theoremstyle{definition}

\theoremstyle{remark}

\def\mm{{\mathfrak m}}

\catcode`\á=\active \def á{\'a}
 \catcode`\Á=\active \def Á{\'A}
  \catcode`\ó=\active \def ó{\'o}
  \catcode`\é=\active \def é{\'e}
  \catcode`\ú=\active \def ú{\'u}
  \catcode`\í=\active \def í{\'{\i}}
  \catcode`\ñ=\active \def ñ{\~n}
  \catcode`\Ñ=\active \def Ñ{\~N}
  \catcode`\¿=\active \def ¿{?`}
  \catcode`\º=\active \def º{$^{\underline{o}}$}
  \catcode`\ª=\active \def ª{$^{\underline{a}}$}
  \catcode`\¡=\active \def ¡{!`}
  \catcode`\â=\active \def â{\^{a}}
  \catcode`\ê=\active \def ê{\^{e}}
  \catcode`\î=\active \def î{\^{\i}}
  \catcode`\ô=\active \def ô{\^{o}}
  \catcode`\û=\active \def û{\^{u}}
  \catcode`\ç=\active \def ç{\c{c}}
  \catcode`\ü=\active \def ü{\"{u}}
  \catcode`\ö=\active \def ö{\"{o}}

 \newenvironment{dem}{\noindent{\it Proof}.}{\qed}

\newcommand{\pro}{\mathbb{P}}
\newcommand{\oo}{\mathcal{O}}
\newcommand{\ooo}{\overline{\mathcal{O}}}
\newcommand{\ago}{\mathfrak{a}}
\newcommand{\bgo}{\mathfrak{b}}
\newcommand{\cgo}{\mathfrak{c}}

\newcommand{\Nu}{{ \mbox{{\LARGE $\nu$}} } }

      \makeatletter
      \def\@setcopyright{}
      \def\serieslogo@{}
      \makeatother

\begin{document}

\author{Julio Jos\'e Moyano-Fern\'andez}
\address{Institut f\"ur Mathematik, Universit\"at Osnabr\"uck. Albrechtstra\ss e 28a, D-49076 Osnabr\"uck, Germany}
\email{jmoyanof@uni-osnabrueck.de}

\title[Fractional ideals and integration]{Fractional ideals and integration with respect to the generalised Euler characteristic}

\begin{abstract}
Let $\bgo$ be a fractional ideal of a one-dimensional
Cohen-Macaulay local ring $\oo$ containing a perfect field $k$.
This paper is devoted to the study some $\oo$-modules associated
with $\bgo$. In addition, different motivic Poincaré series are
introduced by considering ideal filtrations associated with
$\bgo$; the corresponding functional equations of these Poincaré
series are also described.
\end{abstract}



 \thanks{The author was partially supported by the Spanish Government Ministerio de Educaci\'on y Ciencia (MEC)
grant MTM2007-64704 in cooperation with the European Union in the framework of the founds ``FEDER'', and by the Deutsche 
Forschungsgemeinschaft (DFG)}

\maketitle

\section{Introduction}

Canonical ideals over one-dimensional Cohen-Macaulay rings were
profusely studied by Kunz, Herzog, et al., giving nice
characterisations of the Gorenstein property: the formula bringing
conductor and delta-invariant together, or relating it with
symmetry properties of the value semigroup of the ring (cf. e.g. \cite{kunz}; \cite{herzogkunz}). This idea
was developed by Delgado for several branches 
of complex curve singularities (see \cite{Delgado1}), and by Campillo, Delgado and Kiyek
in a more general context (see \cite{cadeki}); they also introduced a
Poincaré series $P(\underline{t})$, deducing its functional
equations when the ring is Gorenstein. Later, Campillo, Delgado
and Gusein-Zade showed that this Poincaré series coincides with
the integral over the projectivisation of the ring with respect to
the Euler characteristic (cf. \cite {cadegu1},\cite{cadegu4}). In
fact, a motivic approach to $P(\underline{t})$ was also introduced
just by taking in the integral the generalised Euler characteristic instead.
Recently, some connections with number-theoretical local
zeta functions were founded by specialising to the case of finite
fields (cf. \cite{demo}). In that paper, the authors introduce
integrals with respect to the generalised Euler characteristic
over fractional ideals of the ring. The aim of the actual work is
to discuss more in detail such motivic integrals, which turn out to be
Poincaré series of fractional ideals, and related objects; in
particular, we extend some of the results concerning the
``classical" series $P(\underline{t})$ and the Gorenstein property
contained in \cite{herzogkunz} and \cite{cadeki}. We will also use
techniques of motivic integration.
\medskip

The paper goes as follows. Section \ref{sec:dos} is an
introduction to the notion of canonical ideal. We define the dual
of an ideal, characterise the self-dual ideals and show that it is indeed
an extension of well-known properties of Gorenstein rings (Theorem
\text{(2.12)}). Section \ref{Sb} is devoted to study an
analogous of the value semigroup $S(\oo)$ of a one-dimensional
Cohen-Macaulay local ring $\oo$ for a fractional ideal
$\mathfrak{b}$: the resulting set $S(\mathfrak{b})$ is no longer a
semigroup, but it has structure of module over $S(\oo)$. We will
give a notion of the symmetry of $S(\mathfrak{b})$; the statements
of Section \ref{sec:dos} allow us to characterise the symmetry of
$S(\mathfrak{b})$ (cf. Proposition \text{(3.8)}). In
Section \ref{generalPS} we define a Poincaré series of motivic
nature for the fractional ideal $\mathfrak{b}$ (Definition
\text{(4.10)}); we prove its rationality and also its functional
equations in absence of the Gorenstein condition (Theorem
\text{(4.19)}). Finally, in Section \ref{functional2} we
investigate the analogous of the extended semigroup
$\widehat{S}(\oo)$ of the ring $\oo$ for a fractional ideal--it
has again structure of $\widehat{S}(\oo)$-module--, as well as an
alternative motivic Poincaré series associated with
$\mathfrak{b}$, giving its functional equations as well (see
Proposition \text{(5.5)}, Theorem \text{(5.8)}).

\section{Duality and fractional ideals} \label{sec:dos}

\textbf{(2.1)} \label{para:21} Let $\oo$ be a one-dimensional Cohen-Macaulay local ring
containing a perfect field $k$ with maximal ideal $\mm$. Let
$\ooo$ be its integral closure with respect to its total ring of
fractions $\mathcal{K}$. Let us assume that $\ooo$ is a finitely
generated $\oo$-module and that the degree $\rho:=[\oo / \mm : k]$
is finite. Let $\delta:= \dim_{k} \left ( \ooo / \oo \right )$ be
the $\delta$-invariant of the ring $\oo$.
\medskip

The ring $\ooo$ decomposes into a finite intersection of Manis
valuation rings, let us say $\overline{\oo} = V_1 \cap \ldots \cap
V_r$. If $\mm(V_i)$ denotes the maximal ideal of $V_i$ for every
$i \in \{1, \ldots , r \}$, then the ideals $\mm_i:= \mm (V_i)
\cap \overline{\oo}$ are principal, regular and maximal (cf.
\cite[Theorem II.(2.11)]{kiyek}) so that $\mm_i=t_i
\overline{\oo}$ for every $i \in \{1, \ldots , r \}$. If $k_i :=
V_i / \mm (V_i) = \overline{\oo} / \mm_i$ for every $i \in \{ 1,
\ldots , r\}$, then the extension degrees
\[
d_i := [k_i : k], \ \ \ \  i \in \{1, \ldots ,r\}
\]
are finite (because $\overline{\oo}$ is a finitely generated
$\oo$-module). Let us also define $d:=d_1 + \ldots + d_r$.

\textbf{(2.2)} \label{rem:cero} Recall that a fractional ideal of a ring $R$
is a $R$-submodule $\mathfrak{a} \ne (0)$ of the total ring of
quotients of $R$ such that $a \mathfrak{a} \subseteq R$ for some
$a \in \oo \setminus \{0\}$. Let us take a fractional ideal
$\mathfrak{a}$ of $\ooo$. It is easy to see that $\mathfrak{a}$
can be written uniquely as a product $\mathfrak{a}=\mm_1^{v_1}
\cdot \ldots \cdot \mm_r^{v_r}$ for $\underline{v}:=(v_1, \ldots ,
v_r) \in \mathds{Z}^r$. We will denote
$\mm^{\underline{v}}:=\mm_1^{v_1} \cdot \ldots \cdot \mm_r^{v_r}$
and $\underline{v}(\mathfrak{a}):=v_1 + \ldots + v_r$.
\medskip

\textbf{(2.3)~Definition:} 
{\sl A fractional ideal $\mathfrak{c}$ of $\oo$ is called
\emph{canonical} if it satisfies the following two
properties:
\begin{itemize}
    \item[(a)] $\mathfrak{c} \cdot \mathcal{K} = \mathcal{K}$
    (i.e., $\mathfrak{c}$ is a regular ideal of $\oo$).
    \item[(b)] For any regular fractional ideal $\ago$ of $\oo$ one
    has that 
    \[
    \ago= \mathfrak{c}:(\mathfrak{c} : \ago).
    \]
\end{itemize}
}

\textbf{(2.4)}  \label{lem:existencia} Since $\oo$ is a one-dimensional local
ring and $\ooo$ is a finitely generated $\oo$-module, a canonical
ideal of $\oo$ does always exist (it is a consequence of
\cite[Satz 2.9, p. 22]{herzogkunz} and \cite[Korollar 2.12, p.
24]{herzogkunz}).

\textbf{(2.5)} Let us fix a canonical ideal $\mathfrak{c}$ of $\oo$. For
every regular fractional ideal $\ago$ of $\oo$, we shall denote
$\ago^{\ast}:=(\mathfrak{c}: \ago)$. The ideal $\ago^{\ast}$ will
be called the dual (ideal) of $\ago$. The ideal $\ago$ is said to
be self-dual if $\ago=\ago^{\ast}$.

\textbf{(2.6)~Lemma:} 
{\sl 
The ring $\oo$ is a canonical ideal of $\oo$ if and only if $\oo$
is self-dual.
}

\dem~ If $\oo$ is canonical, it belongs to the same class (modulo
linear equivalence) as the ideal $\mathfrak{c}$. Then
$\oo^{\ast}=\mathfrak{c}:\oo = \oo:\oo=\oo$. Conversely, assume
$\oo=\mathfrak{c}:\oo$; we have
$\mathfrak{c}:\mathfrak{c}=\mathfrak{c}^{\ast}=\oo$ (\cite[Bemerkung 2.5, p.19]{herzogkunz}), therefore
$\mathfrak{c}=\oo^{\ast}$. \qed

\textbf{(2.7)} \label{label:gore} Notice that the conditions of \text{(2.6)} are
equivalent to the Gorenstein condition (cf. \cite[Korollar 3.4, p.
27]{herzogkunz}). There is also a useful numerical equivalence
given by Theorem \text{(2.9)}.

\textbf{(2.8)} Let us denote by $\lambda_{\oo} (\cdot) = \lambda (\cdot)$ the
length of a finite $\oo$-module. Moreover, we denote
$\mathfrak{f}:=(\oo:\ooo)$; it will be called the conductor ideal
of $\oo$ in $\ooo$, and we can easily check that it is the biggest
ideal of $\oo$ and $\ooo$ at the same time. Let us define $\underline{v}(\mathfrak{f})=:\gamma$ and
\[
\gamma^{\bgo}:= \underline{v} \left ( (\bgo : \ooo): \bgo \right )
= \underline{v} (\bgo : \ooo) - \underline{v} (\bgo \cdot
    \ooo).
\]
Notice that $\gamma^{\oo}=\gamma$.

\textbf{(2.9)~Theorem:~(Gorenstein;~Ap\'ery;~Samuel;~Herzog,~Kunz)} 
{\sl 
We have:
\[
2 \lambda (\oo / \mathfrak{f}) \le \lambda (\ooo /
\mathfrak{f}).
\]
Moreover, the equality holds if and only if the ring $\oo$ is Gorenstein.
}

\textbf{(2.10)} The rest of the section is devoted to generalise Theorem
\text{(2.9)} to any fractional ideal of $\oo$. The obvious task
is to find a candidate to substitute the conductor ideal in the
formula preserving such dimensions. It is easy now: the ideal
$\bgo : \ooo$ is the biggest fractional $\ooo$-ideal in
$\mathfrak{b}^{\ast}$ and the ideal $\bgo \cdot \ooo$ is the
smallest $\ooo$-ideal containing $\bgo$ and $\ooo$. Notice also
that $(\bgo:\ooo)^{\ast}=\bgo^{\ast} \cdot \ooo$ and $(\bgo \cdot
\ooo)^{\ast}=(\bgo^{\ast}:\ooo)$. First of all note the following fact (cf. \cite{stohr}):

\textbf{(2.11)~Lemma:} 
{\sl 
We have
\[
2 \lambda (\mathfrak{b}^{\ast} \cdot \oo / \mathfrak{b}^{\ast} ) = 2 \lambda (\mathfrak{b} / \mathfrak{b}:\ooo ).
\]
}

\textbf{(2.12)~Theorem:} 
{\sl 
Let $\mathfrak{b}$ be a fractional ideal of $\oo$. We have
\[
2 \lambda (\mathfrak{b} / \mathfrak{b}:\ooo ) \le \lambda
(\mathfrak{b} \cdot \ooo /\mathfrak{b}:\ooo).
\]
Moreover, the equality holds if and only if $\mathfrak{b}$ is
self-dual.
}

\dem~Without loss of generality, let us assume that $\bgo
\subseteq \oo \subseteq \mathfrak{c} \subseteq \ooo$. Now
$\bgo^{\ast} = \cgo:\bgo \supseteq \bgo$. Indeed $\bgo^{\ast}
\supseteq \cgo^{\ast}=\oo \supseteq \bgo$. Then we have
\[
\xymatrix{\bgo : \ooo~ \ar@{^{(}->}[r] & \bgo~ \ar@{^{(}->}[r] \ar@{^{(}->}[dr] & \oo~ \ar@{^{(}->}[r] &  \bgo^{\ast}~ \ar@{^{(}->}[r] & \bgo^{\ast}\cdot \ooo   \\
                          &                  & \bgo \cdot \ooo \ar@{^{(}->}"2,3":"1,5"             &                      &
                          }\eqno(\dagger)
\]
By looking at
($\dagger$) we have
\begin{align*}
\lambda (\bgo^{\ast} \cdot \ooo/\bgo : \ooo) = & \lambda
(\bgo^{\ast} \cdot \ooo / \bgo \cdot \ooo) + \lambda (\bgo \cdot
\ooo / \bgo : \ooo)
\\
= & \lambda (\bgo^{\ast} \cdot \ooo / \bgo \cdot \ooo ) + \lambda
(\bgo \cdot \ooo / \bgo) + \lambda (\bgo / \bgo : \ooo).
\end{align*}

If $\bgo$ is self--dual, then $\lambda (\bgo \cdot \ooo / \bgo) = \lambda (\bgo / \bgo : \ooo)$ by Lemma \text{(2.11)}, and by substituting above we are done. Conversely, assuming the following equalities hold:
\[
 \lambda
(\mathfrak{b} \cdot \ooo /\mathfrak{b}:\ooo) \overset{(1)}{=} 2\lambda (\mathfrak{b} / \mathfrak{b}:\ooo ) \overset{(2)}{=} 2\lambda (\mathfrak{b}^{\ast} \cdot \oo / \mathfrak{b}^{\ast});
\]
again looking at ($\dag$) we get
\[
 \lambda (\bgo^{\ast} \cdot \ooo /\bgo \cdot \ooo) \overset{(3)}{=} \lambda
(\bgo^{\ast} \cdot \ooo / \bgo^{\ast}) + \lambda (\bgo^{\ast}/\bgo) + \lambda (\bgo/\bgo : \ooo).
\]
By (1) and (2) 
\[
 \lambda (\bgo^{\ast} \cdot \ooo /\bgo \cdot \ooo) = \lambda (\mathfrak{b} / \mathfrak{b}:\ooo ) + \lambda (\mathfrak{b} / \mathfrak{b}:\ooo ) = \lambda (\mathfrak{b} / \mathfrak{b}:\ooo ) + \lambda (\mathfrak{b}^{\ast} \cdot \oo / \mathfrak{b}^{\ast}),
\]
and plugging the last equality into (3) shows $\lambda (\bgo^{\ast}/\bgo)=0$, i.e., $\bgo^{\ast}=\bgo$, hence the ideal $\bgo$ is self--dual. \qed


%
\textbf{(2.13)~Remark:} 
{\sl
From Lemma \text{(2.11)} and Theorem \text{(2.12)} it follows: The ideal $\bgo$ is self--dual if and only if
\[
2 \lambda (\mathfrak{b} \cdot \oo / \mathfrak{b}) = 2 \lambda (\mathfrak{b} / \mathfrak{b}:\ooo ) = \lambda
(\mathfrak{b} \cdot \ooo /\mathfrak{b}:\ooo).
\]
}
\textbf{\!\!(2.14)}  Let $\underline{v}:=(v_1, \ldots , v_r), \underline{w}:=(w_1,
\ldots , w_r)$ be vectors in $\mathds{Z}^r$. We will write
$\underline{v} \ge \underline{w}$ if and only if $v_i \ge w_i$ for
every $i \in \{1, \ldots ,r \}$, $\underline{0}:=(0,0, \ldots, 0)
\in \mathds{Z}^r$ and $\underline{1}:=(1,1, \ldots ,1) \in
\mathds{Z}^r$. Moreover, for every subset $I \subseteq I_0:=\{1,
\ldots , r\}$ let $\sharp I$ be the number of elements in $I$, and
let $\underline{1}_{I}$ be the element of $\mathds{Z}^r$ whose
$i$th component is equal to $1$ or $0$ if $i \in I$ or $i \notin
I$ respectively. For any $\underline{v} \in \mathds{Z}^r$ and any
fractional ideal $\bgo$ in $\oo$, we define the set
\[
J^{\bgo}(\underline{v}):=\{z \in \bgo \setminus \{0\} \mid
\underline{v}(z) \ge \underline{v}\},
\]
with $\underline{v}(z):=(v_1(z), \ldots, v_r(z))$. They are ideals
defining a multi--index
filtration $\{J^{\bgo}(\underline{v}) \}$, as
$J^{\bgo}(\underline{v}) \supseteq J^{\bgo}(\underline{w})$ if
$\underline{w} \ge \underline{v}$. 
\medskip

For every $i \in \{1, \ldots ,r\}$, let us define
$C^{\bgo}(\underline{v},i):=J^{\bgo}(\underline{v})/J^{\bgo}(\underline{v}+\underline{1}_{\{i\}})$, and $c^{\bgo}(\underline{v},i):=\dim_k (C^{\bgo}(\underline{v},i))$; we write also 
$C^{\bgo}(\underline{v}):=J^{\bgo}(\underline{v})/J^{\bgo}(\underline{v}+\underline{1})$, as well as
$c^{\bgo}(\underline{v}):=\dim_k
(C^{\bgo}(\underline{v}))$.
Since $\oo$ is Cohen-Macaulay, $c^{\bgo} (\underline{v}) < \infty$
for every $\underline{v} \in \mathds{Z}^r$; also the filtration is
finitely determined, i.e., for any $\underline{v} \in
\mathds{Z}^r$ there
 exists $N \in \mathds{Z}$ such that $J^{\bgo}(\underline{v}) \supset
 \mm^N$; that means that every subspace $J^{\bgo}(\underline{v})$ of
 $\oo$ has finite codimension $\ell^{\bgo} (\underline{v})$ (cf. \cite[p. 194]{cadegu11}). Notice that, for every $i \in \{1, \ldots , r\}$ one has $0 \le c^{\bgo}(\underline{v},i) \le d_i$, and if $\underline{v} \ge \gamma^{\bgo}$, then $c^{\bgo}(\underline{v},i) = d_i$ for every $i \in \{1, \ldots , r\}$ (The proof of these facts follows much more \cite{cadeki}).

\textbf{(2.15)} 
The Gorenstein condition on the ring $\oo$ was proven to be equivalent
 to the following equality (cf. \cite[Corollary (3.7)]{cadeki}):
 \[
c^{\oo} (\underline{v}) +
c^{\oo}(\gamma-\underline{v}-\underline{1}) =d, ~ ~ ~
\mathrm{~for~every~} \underline{v} \in \mathds{Z}^r.
 \]
Our purpose now is to state the analogue of this result for the
case of a fractional ideal. The proof is adapted from
\cite{cadeki}. We show first:

\textbf{(2.16)~Lemma:} 
{\sl 
We have $c^{\bgo}(\gamma^{\bgo}-\underline{1}_{\{i\}},i) < d_i$.
}

\dem~ Write $\gamma^{\bgo}=(\gamma_1^{\bgo}, \ldots , \gamma_r^{\bgo})$. Let $i \in \{1,
\ldots ,r\}$. Consider the $k$-linear map $\phi_i: C^{\bgo}(\gamma^{\bgo}-\underline{1}_{\{i\}})  \to \bgo \cdot \ooo / \bgo : \ooo$ given by 
\[
z \mathrm{~mod~} J^{\bgo}(\gamma^{\bgo}-\underline{1}_{\{i\}}) \mapsto z \cdot t_i^{-(\gamma_i^{\bgo}-1)} \mathrm{~mod~} \bgo : \ooo. 
\]
This map is clearly injective, so we have to prove that $\phi$ is not an epimorphism. Let $i \in \{1, \ldots , r\}$. It is easily seen that 
\[
t^{\gamma^{\bgo}-\underline{1}_{\{i\}}} \notin (\bgo : \ooo):(\bgo \cdot \ooo)=(\bgo : \ooo): \bgo.
\]
Hence there exists $\xi \in \bgo \cdot \ooo$ such that $\xi \cdot t^{\gamma^{\bgo}-\underline{1}_{\{i\}}} \notin \bgo$. Let $x_i:=\xi \cdot t^{\gamma^{\bgo}}\cdot t_i^{-\gamma_i^{\bgo}}$. Notice that $x_i \in \bgo \cdot \ooo$. Furthermore, if $z \in J^{\bgo} (\gamma^{\bgo}-\underline{1}_{\{i\}})$ then $\xi \cdot t^{\gamma^{\bgo}-\underline{1}_{\{i\}}} - z \notin \bgo : \ooo = \bgo \cdot \mm^{\gamma^{\bgo}}$ by definition of $J^{\bgo} (\gamma^{\bgo}-\underline{1}_{\{i\}})$; thus $v_i (\xi \cdot  t^{\gamma^{\bgo}-\underline{1}_{\{i\}}}-z)=\gamma_i^{\bgo}-1+v_i(\bgo \cdot \ooo)$ and so 
\[
v_i(x_i-\xi \cdot t_i^{-(\gamma_i^{\bgo}-1)})=v_i (\xi \cdot t^{\gamma^{\bgo}-\underline{1}_{\{i\}}} \cdot t_i^{-\gamma_i^{\bgo}}-z \cdot t_i^{1-\gamma_i^{\bgo}})=v_i(\bgo \cdot \ooo), 
\]
which proves the non--surjectivity of $\phi_i$.

\textbf{(2.17)~Proposition:} 
{\sl 
Assume $\oo$ to be analytically reduced. Let $\bgo$ be a
fractional ideal, let $\underline{v} \in \mathds{Z}^r$; then
\[
c^{\bgo}(\underline{v},i) + c^{\bgo}(\gamma^{\bgo} - \underline{v} -\underline{1}_{\{i\}},i)
\le d_i ~ ~ ~ ~ \mathrm{~for~every~} i \in \{1, \ldots, r\}.
\]
}

\dem~
For every $\underline{v} \in \mathds{Z}^r$ and every $i \in \{1,
\ldots ,r\}$, the map
$\eta_{\underline{v},i}:C^{\bgo}(\underline{v},i) \to k_i$ defined
by $z \mathrm{~mod~} J^{\bgo}(\underline{v}+\underline{1}_{\{i\}}) \mapsto z t_i^{-v_i} \mathrm{~mod~}
\mm_i$ is a $k$-monomorphism, hence
$c^{\bgo}(\underline{v},i) \le d_i$. Take now
$\underline{w}=(w_1, \ldots ,w_r) \in \mathds{Z}^r$ with
$\underline{v} \le \underline{w}$. Thus for every $i \in \{1,
\ldots, r\}$, the inclusion $J^{\bgo}(\underline{w}) \to
J^{\bgo}(\underline{v})$ induces a $k$-homomorphism
$\varphi_{\underline{w},\underline{v},i}:C^{\bgo}(\underline{w},i)
\to C^{\bgo}(\underline{v},i)$. For an $i \in \{1, \ldots, r\}$
with $v_i=w_i$ we have
$\eta_{\underline{w},i}=\eta_{\underline{v},i} \circ
\varphi_{\underline{w},\underline{v},i}$ and
$\varphi_{\underline{w},\underline{v},i}$ is injective. Notice
also that for every $\underline{n},\underline{m} \in \mathds{Z}^r$
the following inclusion holds:
\[
J^{\bgo}(\underline{v})J^{\oo}(\underline{w})+J^{\bgo}(\underline{w})J^{\oo}(\underline{v})
\subseteq J^{\bgo}(\underline{v}+\underline{w}).\eqno(\ast)
\]
Then for $a \in J^{\bgo}(\underline{v}), b \in
J^{\bgo}(\underline{w}), i \in \{1, \ldots ,r\}$ we have
\[
\eta_{\underline{v},i}(a \mathrm{~mod~}
J^{\bgo}(\underline{v}+\underline{1}_{\{i\}})) \cdot
\eta_{\underline{w},i}(b
\mathrm{~mod~}J^{\bgo}(\underline{w}+\underline{1}_{\{i\}})) =
\eta_{\underline{v}+\underline{w},i}(J^{\bgo}(ab
\mathrm{~mod~}\underline{v}+\underline{w})).
\]
Let $i \in \{1, \ldots ,r\}$. Let $H_i \subset k_i$ be a
$1$-codimensional subspace of the $k$-vector space $k_i$
containing
$\mathrm{im}(\eta_{\gamma^{\bgo}-\underline{1}_{\{i\}},i})$ (notice that it is possible by Lemma \text{(2.16)}); consider the
$k$-bilinear pairing $k_i \times k_i \to k_i \to k_i/H_i$ defined
by $(a,b) \to a \cdot b \mathrm{~mod~} H_i$ (cf.
\cite[(3.5)]{cadeki}), which is non--degenerate (multiplying by scalars of $k_i$ is a $k$-automorphism on $k_i$). Because of $(\ast)$ we have
\[
J^{\bgo}(\underline{v})J^{\oo}(\gamma^{\bgo}-\underline{v}-\underline{1}_{\{i\}})
+ J^{\bgo}(\gamma^{\bgo}-\underline{v}-\underline{1}_{\{i\}})
J^{\oo}(\underline{v}) \subset
J^{\bgo}(\gamma^{\bgo}-\underline{1}_{\{i\}}),
\]
therefore $\mathrm{im}(\eta_{\underline{v},i})$ lies in the
orthogonal complement of
$\mathrm{im}(\eta_{\gamma^{\bgo}-\underline{v}-\underline{1}_{\{i\}},i})$,
hence
\[
c^{\bgo}(\underline{v},i) \le d_i - c^{\bgo}(\gamma^{\bgo}-\underline{v}-\underline{1}_{\{i\}},i).
\]
\qed

It remains to show that the equality of Proposition \text{(2.17)} holds if and only if $\bgo$ is
self-dual. This follows by the same method as in the proof of
\cite[Theorem (3.6)]{cadeki}, just by applying the
characterisation of the self-dual fractional ideals provided by
Theorem \text{(2.12)} instead of using Theorem
\text{(2.9)}. 

\textbf{(2.18)~Lemma:} 
{\sl 
Let $\bgo \cdot \ooo = \ooo$. Let $\{v^{(p)}\}_{0 \le p \le h}$ be a strictly increasing sequence in $\mathds{Z}^r$ such that $v^{(0)}=0$, $v^{(h)}=\gamma^{\bgo}$, and for every $p \in \{1, \ldots , h\}$ there exists $i(p) \in \{1, \ldots , r\}$ satisfying $v^{(p)}-v^{(p-1)} = \underline{1}_{\{i(p)\}}$. Then $\bgo$ is self--dual if and only if 
\[
c^{\bgo}(v^{(p)},i(p+1))+c^{\bgo}(\gamma^{\bgo}-v^{(p)}-\underline{1}_{\{i(p+1)\}},i(p+1))=d_{i(p+1)}
\]
for every $p \in \{0, \ldots , h-1\}$.
}

\dem~ Define $w^{(p)}:=\gamma^{\beta}-v^{(p)}$ for every  $p \in \{0, \ldots , h\}$. We have $v^{(p)}+w^{(p+1)}=\gamma^{\bgo}-\underline{1}_{\{i(p+1)\}}$, $w^{(p)}+w^{(p+1)}=\underline{1}_{\{i(p+1)\}}$ and
\[
\bgo = J^{\bgo}(v^{(0)}) \supset J^{\bgo}(v^{(1)}) \supset \ldots \supset J^{\bgo}(v^{(h)})= \bgo :\ooo.
\]
\[
\bgo : \ooo= J^{\bgo}(w^{(0)}) \subset J^{\bgo}(w^{(1)}) \subset \ldots \subset J^{\bgo}(w^{(h)})= \bgo.
\]
Therefore
\[
\sum_{p=0}^{h-1} c^{\bgo}(v^{(p)},i(p+1))=\dim_k (\bgo / \bgo: \ooo) = \sum_{p=0}^{h-1} c^{\bgo}(w^{(p)},i(p+1)),
\]
and then
\[
2 \dim_k (\bgo / \bgo:\ooo)=\sum_{p=0}^{h-1} c^{\bgo}(v^{(p)},i(p+1))+c^{\bgo}(w^{(p)}-\underline{1}_{\{i(p+1)\}},i(p+1)).
\]
By Proposition \text{(2.17)} this expression is smaller than or equal to
\[
\sum_{p=0}^{h-1} d_{i(p+1)}=\sum_{p=1}^{r} \gamma_i^{\bgo} d_i = \dim_k (\ooo / \bgo : \ooo),
\]
and by Theorem \text{(2.12)}, and the assumption $\bgo \cdot \ooo=\ooo$, we are done. The converse follows in the same manner as in the proof of \cite[(3.6)]{cadeki}, part (d).
\qed

\textbf{(2.19)~Theorem:} 
{\sl 
Let $\bgo$ be a fractional $\oo$-ideal such that $\bgo \cdot \ooo = \ooo$.
The following statements are equivalent:
\begin{itemize}
\item[(1)] For every $\underline{v} \in \mathds{Z}^r$ we have $ c^{\bgo}(\underline{v},i)+c^{\bgo}(\gamma^{\bgo}-\underline{v}-\underline{1}_{\{i\}},i) \overset{\ast}{=}d_i$,
for $i \in \{1, \ldots , r\}$.
\item[(2)] $\bgo$ is self--dual.
\end{itemize}
}

\dem~ If the equality $\ast$ holds for every $\underline{v} \in \mathds{Z}^r$ and for every $i \in \{1, \ldots , r\}$, then one can choose a strictly increasing sequence as in Lemma \text{(2.18)}, and we obtain that $2 \dim (\bgo / \bgo:\ooo)=\dim (\ooo/\bgo: \ooo)$, which by Theorem \text{(2.17)}  implies the statement. \qed

Analogously as in \cite[(3.7)]{cadeki} one shows:

\textbf{(2.20)~Corollary:} 
{\sl 
Let $\underline{v} \in \mathds{Z}^r$. Then 
\[
c^{\bgo}(\underline{v})+c^{\bgo} (\gamma^{\bgo}-\underline{v}-\underline{1}) \le d=\sum_{i=1}^{r}d_i.
\]
Moreover, $\bgo$ is self--dual if and only if the equality holds for every $\underline{v} \in \mathds{Z}^r$.
}

\section{The value ideal of a fractional ideal} \label{Sb}

Let $\bgo$ be a fractional ideal of $\oo$. Consider the set
\[
S(\bgo)= \{\underline{v}(z)-\underline{v}(\bgo \cdot \ooo) \mid z \in \bgo \setminus \{0\}\}.
\]

This is a subsemigroup of $\mathds{Z}^r$ which is in fact a
$S(\oo)$-module and only depends on the ideal class of $\bgo$ in the ideal class semigroup of $\oo$. The elements of $S(\bgo)$ are connected to
filtrations $\{J^{\bgo}(\underline{v})\}$ in the following sense:

\textbf{(3.1)~Lemma:} 
{\sl 
Let $k$ be an infinite field. Let $\underline{v}=(v_1, \ldots
,v_r) \in \mathds{Z}^r$. Then $\underline{v} \in S(\bgo)$ if and
only if
$J^{\bgo}(\underline{v})/J^{\bgo}(\underline{v}+\underline{1}_{\{i\}})
\ne 0$ for every $i \in \{1, \ldots,r\}$.
}

\dem~ If $\underline{v} \in S(\bgo)$, then there exists
$\underline{w}=(w_1, \ldots,w_r) \in S(\bgo)$ such that $w_i=v_i$
and $w_j \ge v_j$ for every $j \in \{1, \ldots ,r\}$, $j \ne i$.
Thus
$J^{\bgo}(\underline{v})/J^{\bgo}(\underline{v}+\underline{1}_{\{i\}})
\ne 0$ for every $i \in \{1, \ldots,r\}$. Conversely, assume that
$J^{\bgo}(\underline{v})/J^{\bgo}(\underline{v}+\underline{1}_{\{i\}})
\ne 0$ for every $i \in \{1, \ldots,r\}$. For every $i \in \{1,
\ldots ,r\}$ choose an element $z_i \in J^{\bgo}(\underline{v})
\setminus J^{\bgo}(\underline{v}+\underline{1}_{\{i\}})$. Since
$k$ is infinite and $k=k_i$ for every $i \in \{1, \ldots ,r\}$,
then there exist elements $a_1, \ldots , a_r \in \oo$ such that
$v_i(a_1z_1+ \ldots + a_rz_r)=v_i $ for every $i \in \{1, \ldots ,
r\}$, i.e., $\underline{v}(a_1z_1+ \ldots + a_r
z_r)=\underline{v}$.\qed

\textbf{(3.2)} Notice that if $\bgo=\oo$ then $S(\oo)$ is the value semigroup
of the ring $\oo$. Recall that $\gamma^{\bgo}:= \underline{v} \left ( (\bgo : \ooo): \bgo \right )
= \underline{v} (\bgo : \ooo) - \underline{v} (\bgo \cdot \ooo)$, and
$\gamma^{\oo}=\gamma$ if $\bgo=\oo$. The next two lemmas
show that $\gamma^{\bgo}$ plays the role of a kind of conductor of
the set $S(\bgo)$.

\textbf{(3.3)~Lemma:} 
{\sl 
For every fractional ideal $\bgo$ in $\oo$, there exists
$\underline{m} \in \mathds{Z}^r$ such that
$J^{\mathcal{K}}(\underline{n}) \subseteq \bgo$ for all
$\underline{n} \ge \underline{m}$.
}

\dem~ First of all, notice that $(\oo:\ooo)= \{ x \in \ooo \mid
\underline{v}(x) \ge \gamma \}$. Then $x \in z^{-1}\oo$ for all $x
\in \ooo$ and $z \in \bgo$, hence $(\oo:\ooo) \subseteq z^{-1}
\bgo$ for all $z \in \bgo$ and $z \cdot (\oo:\ooo) = \{x \mid
\underline{v}(x) \ge \gamma + \underline{v}(z) \} \subseteq \bgo$
for all $z \in \bgo$. Therefore $J^{\mathcal{K}}(\underline{n})
\subset \bgo$ for every $\underline{n} \ge \underline{v}(z)$.\qed


\textbf{(3.4)~Lemma:} 
{\sl 
We have:
\begin{itemize}
    \item[(a)] $\underline{0} \le \gamma^{\bgo} \le \underline{v} \left ( (\oo:\ooo) \right ) =
    \gamma$.
    \item[(b)] If $\bgo \cdot \ooo=\ooo$, then $\gamma^{\bgo}:= \min \{ \underline{n} \mid \mm^{\underline{n}}
\subseteq \bgo \}$.
\end{itemize}
}

\dem~(a)~From \cite[Lemma 3.1]{stohr3}, one has the inclusions
\[
(\oo:\ooo) \subseteq (\bgo : \ooo):(\bgo \cdot \ooo) \subseteq \ooo
\]
and the statement follows.~ The assertion (b)~can be deduced from the fact that
$(\bgo:\ooo):\bgo=(\bgo:\ooo):(\bgo \cdot \ooo)$. \qed

\textbf{(3.5)} In the rest of the section, the ring $\oo$ will be assumed to
be residually rational, i.e., $k=k_i$ for every $i \in \{1, \ldots
,r\}$.

\textbf{(3.6)} For every $\underline{n}=(n_1, \ldots , n_r) \in \mathds{Z}^r$
and for every $i \in \{1, \ldots , r\}$, define
\[
\Delta_i(\underline{n})= \{ (\sigma_1, \ldots , \sigma_r) \in
S(\bgo) \mid \sigma_i=n_i \mathrm{~and~}\sigma_j > \sigma_i
\mathrm{~for~} j \in \{1, \ldots , r\}, j \ne i \}.
\]
Moreover, for every $\underline{n} \in \mathds{Z}^r$ define
\[
\Delta(\underline{n}) = \bigcup_{i=1}^{r} \Delta_i(\underline{n}).
\]

\textbf{(3.7)~Definition:} 
{\sl 
The $S(\oo)$-module $S(\bgo)$ is said to be symmetric if there
exists $\tau \in \mathds{Z}^r$ such that, for every $\underline{v}
\in \mathds{Z}^r$, $\underline{v} \in S(\bgo)$ if and only if
$\Delta(\tau-\underline{v}) = \varnothing$.
}

\medskip

\textbf{(3.8)~Proposition:} 
{\sl 
Let $\oo$ be residually rational. We have
\begin{enumerate}
    \item $\Delta (\gamma^{\bgo}-\underline{v}-\underline{1}) =
    \varnothing$ for every $\underline{v} \in S(\bgo)$.
    \item If $S(\bgo)$ is symmetric then
$\bgo$ is self-dual.
    \item Suppose, in addition, that $k$ is an infinite field; if
    $\bgo$ is a self-dual fractional ideal, then $S(\bgo)$ is symmetric.
\end{enumerate}
}

\dem~(1) Let $\underline{v} \in S(\bgo)$. Then we have
$J^{\bgo}(\underline{v})/J^{\bgo}(\underline{v}+\underline{1}_{\{i\}})
\ne 0$ for every $i \in \{1, \ldots ,r\}$, therefore
$J^{\bgo}(\gamma^{\bgo}-\underline{v}-\underline{1}_{\{i\}})/J^{\bgo}(\gamma^{\bgo}-\underline{v})
\ne 0$ for every $i \in \{1, \ldots ,r\}$ and so
$\Delta(\gamma^{\bgo}-\underline{v}-\underline{1})=\varnothing$.~(2)~Let
$\underline{v}=(v_1, \ldots, v_r) \in \mathds{Z}^r$ and take $i \in
\{1, \ldots ,r\}$. If
$J^{\bgo}(\underline{v})/J^{\bgo}(\underline{v}+\underline{1}_{\{i\}})
\ne 0$ then we must have that
$J^{\bgo}(\gamma^{\bgo}-\underline{v}-\underline{1}_{\{i\}})/J^{\bgo}(\gamma^{\bgo}-\underline{v})
= 0$. Let us consider the case in which 
$J^{\bgo}(\underline{v})/J^{\bgo}(\underline{v}+\underline{1}_{\{i\}})
= 0$. There exists a vector $\underline{w}=(w_1, \ldots,w_r)\in
\mathds{Z}^r$ with $\Delta(\underline{w}) = \varnothing$, $w_i=v_i$
and $w_j < v_j$ for every $j \in \{1, \ldots,r\}$, $j \ne i$.
Since $S(\bgo)$ is symmetric,
$\gamma^{\bgo}-\underline{w}-\underline{1} \in S(\bgo)$. Now
$\gamma^{\bgo}_j-w_j-1 \ge \gamma^{\bgo}-v_j$ for every $j \in
\{1, \ldots ,r\}$, $j \ne i$, and
$\gamma^{\bgo}_i-w_i-1=\gamma^{\bgo}_i-v_i-1$, hence for any
regular element $z \in \oo$ satisfying
$\underline{v}(z)=\gamma^{\bgo}-\underline{w}-\underline{1}$ it
follows that $z \in
J^{\bgo}(\gamma^{\bgo}-\underline{v}-\underline{1}_{\{i\}})$, $z
\notin J^{\bgo}(\gamma^{\bgo}-\underline{v})$, and therefore
$J^{\bgo}(\gamma^{\bgo}-\underline{v}-\underline{1}_{\{i\}})/J^{\bgo}(\gamma^{\bgo}-\underline{v})
\ne 0$. Theorem \text{(2.19)} implies that $\bgo$ is
self-dual.~ (3)~ Assume that $k$ is infinite. Let $\underline{v}
\in \mathds{Z}^r$ and assume that $\Delta
(\gamma^{\bgo}-\underline{v}-\underline{1})=\varnothing$. Then
$J^{\bgo}(\gamma^{\bgo}-\underline{v}-\underline{1}_{\{i\}})/J^{\bgo}(\gamma^{\bgo}-\underline{v})
= 0$ for every $i \in \{1, \ldots ,r\}$, hence
$J^{\bgo}(\underline{v})/J^{\bgo}(\underline{v}+\underline{1}_{\{i\}})$
for every $i \in \{1, \ldots ,r \}$ by Theorem
\text{(2.17)} and therefore $z \in S(\bgo)$ by Lemma
\text{(3.1)}. Thus $S(\bgo)$ is symmetric.\qed

\textbf{(3.9)~Corollary:} 
{\sl 
Let $k$ be an infinite field. The value semigroup $S(\oo)$ is
symmetric if and only if the ring $\oo$ is Gorenstein.
}

\section{Generalised Poincar\'e series of a fractional ideal}
\label{generalPS}

\textbf{(4.1)} Let $\bgo$ be a fractional ideal in $\oo$. The multi-index
filtration $\{J^{\bgo}(\underline{v})\}$ defines a Laurent series
\[
 L(\bgo, t_1, \ldots , t_r):=\sum_{\underline{v} \in \mathds{Z}^r}
 \dim_{k} \left ( J^{\bgo}(\underline{v}) / J^{\bgo}(\underline{v} + \underline{1})\right ) \cdot
 \underline{t}^{\underline{v}} \in \mathds{Z} [\![t_1, \ldots, t_r]\!],
\]
where $\underline{t}^{\underline{v}}:=t_1^{v_1} \cdot \ldots \cdot
t_r^{v_r}$. We will write $L(\bgo,\underline{t})$ instead of
$L(\bgo,t_1, \ldots , t_r)$ if the number of variables is clear
from the context.
\medskip

\textbf{(4.2)} \label{tiposfiltrac} There is a priori no fixed way to choose
a suitable coefficient in $L(\bgo,\underline{t})$. We may consider
the following spaces:
\begin{itemize}
\item[(1)]
$J^{\bgo}(\underline{v})/J^{\bgo}(\underline{v}+\underline{1})$;
\item[(2)]
$J^{\bgo}(\underline{v})/J^{\bgo}(\underline{v}+\underline{1})
\setminus \bigcup_{i=1}^{r}
J^{\bgo}(\underline{v}+\underline{1}_{\{ i
\}})/J^{\bgo}(\underline{v}+\underline{1})$; \item[(3)]
$J^{\bgo}(\underline{v}) \setminus \bigcup_{i=1}^{r}
J^{\bgo}(\underline{v}+\underline{1}_{\{ i \}})$.
\end{itemize}

Filtration $(1)$ is related to the semigroup of values of the
ring, $(2)$ defines the Poincar\'e series in terms of the extended
semigroup of the ring, and $(3)$ introduces the Poincar\'e series
as an integral with respect to the Euler characteristic. Exactly
this last point of view makes clear the association between the
dimension of a vector space and the Euler characteristic $\chi$ of
its projectivisation, namely:
\[
\dim_{k} \left ( J(\underline{v})/J(\underline{v}+\underline{1})
\right ) = \chi \left ( \mathbb{P} \left (
J(\underline{v})/J(\underline{v}+\underline{1}) \right ) \right ).
\]

\textbf{(4.3)} We can also choose other measures than $\chi$, for instance
the so-called generalised Euler characteristic $\chi_g$. It is a
sort of motivic Euler characteristic which makes use of the notion
of Grothendieck ring. The Grothendieck ring $K_0 (\Nu_{k})$ is
defined to be the free Abelian group on isomorphism classes $[X]$
of quasi-projective schemes $X$ of finite type over $k$ subject to
the following relations:
\begin{itemize}
\item[(1)] $[X_1]=[X_2]$ if $X_1 \cong X_2$ for $X_1,X_2 \in
\Nu_{k}$; \item[(2)] $[X]=[X \setminus Z]+[Z]$ for a closed
subscheme $Z$ of $X \in \Nu_k$;
\end{itemize}
and taking the fibred product as multiplication:
\begin{itemize}
\item[(3)] $[X_1] \cdot [X_2]=[X_1 \times_k X_2]$ for $X_1,X_2 \in
\Nu_k$.
\end{itemize}

\textbf{(4.4)} Let $k [T]$ be the polynomial ring in one indeterminate $T$
over the field $k$. The affine scheme $\mathrm{Spec}(k[T])$ over
$k$ is the affine line over $k$, which will be denoted by
$\mathbb{A}^1_k$. The class of the affine line in $K_0(\Nu_k)$, denoted by $\mathbb{L}$, is
called the Lefschetz class of $K_0 (\Nu_k)$.

\textbf{(4.5)} Let $p$ be a non-negative integer and let $J_{\oo}^p$ be the
space of $p$-jets over $\oo$, which is a finite-dimensional
$k$-vector space of dimension $d(p)$. Let us consider its
projectivisation $\pro J_{\oo}^p$ and let us adjoin one point to
this (that is, $\pro^{\ast}J_{\oo}^p = \pro J_{\oo}^p \cup \{\ast
\}$ with $\ast$ representing the added point) in order to have a
well-defined map $\pi_p: \pro \oo \to \pro^{\ast} J_{\oo}^p$. A
subset $X \subset \pro \oo$ is said to be {\sl cylindric} if there
exists a constructible subset $Y \subset \pro J_{\oo}^p \subset
\pro^{\ast} J_{\oo}^p$ such that $X = \pi_p^{-1} (Y)$.

\textbf{(4.6)} The generalised Euler characteristic $\chi_g (X)$ of a
cylindric subset $X$ is the element $[Y] \cdot \mathbb{L}^{-d(p)}$
in the ring $K_0 (\nu_{k})_{(\mathbb{L})}$, where
$Y=\pi^{-1}_p(X)$ is a constructible subset of $\mathbb{P}\oo$.
Note that $\chi_g (X)$ is well-defined, because if
$X=\pi^{-1}_q(Y^{\prime})$, $Y^{\prime} \subset \pro J_{\oo}^q$
and $p \ge q$, then $Y$ is a locally trivial fibration over
$Y^{\prime}$ and therefore $[Y] = [Y^{\prime}] \cdot
\mathbb{L}^{d(p)-d(q)}$.

\textbf{(4.7)} As in \cite{demo}, we can extend these definitions to subsets
of $\mathcal{K}$ (in particular to fractional ideals): a subset $X
\subseteq \mathcal{K}$ is called cylindric if there exists a
non-zero divisor element $z \in \oo$ such that the set $zX$ is a
subset of $\oo$ and is cylindric. In this situation, the
generalised Euler characteristic is
\begin{equation} \label{equation:chiX}
\chi_g(X):= \frac{\chi_g (zX)}{\chi_g (z \oo)}. \nonumber
\end{equation}

Let $\ago \subseteq \oo$ be an ideal of $\oo$. Since $\ago$ is
$\mm$-primary, we have $\mm^{p+1} \subseteq \ago$. Let
$\overline{\ago}$ be the ideal $\ago / \mm^{p+1}$ of $\oo /
\mm^{p+1}$ so that $\pi_p^{-1} (\overline{\ago})=\ago$. As $\oo /
\mm^{p+1}$ is a finite-dimensional $k$-vector space, the ideal
$\overline{\ago}$ is constructible. Then $\ago$ is cylindric and
we get
\begin{eqnarray}
\chi_g (\ago) & = & \left [ \overline{\ago} \right ] \cdot
\mathbb{L}^{-d(p)} \nonumber \\
& = & \mathbb{L}^{\dim_k \left ( \ago / \mm^{p+1} \right )-d(p)}
\nonumber \\
& = & \mathbb{L}^{\deg (\ago)}. \nonumber
\end{eqnarray}
In particular, $\chi_g (\mm^{p+1}) = \mathbb{L}^{-d(p)}$.

\textbf{(4.8)} Let $G$ be an abelian group with countable many values. Let
$X$ be a cylindric subset of $\mathcal{K}$. A function $\psi: X
\to G$ is called cylindric if the set $\psi^{-1}(a) \subseteq
\mathcal{K}$ is cylindric for all $a \in G \setminus \{ 0\}$. The
integral of $\psi$ over $X$ with respect to the generalised Euler
characteristic is
\[
\int_X \psi d \chi_g := \sum_{a \in G \setminus \{0\}} \chi_g
(\psi^{-1}(a)) \cdot a,
\]
if this sum makes sense in $K_0 (\Nu_{k})_{(\mathbb{L})}
\otimes_{\mathds{Z}} G$; in such a case, the function $\psi$ is
said to be \emph{integrable}.
\medskip

\textbf{(4.9)~Remark:} 
{\sl 
Let $\psi: \mathbb{P} X \to G$ be a cylindric function of $X$. Let
us denote by $\psi^{\prime}: X \to G$ the function induced by
$\psi$ on $X$ with $\psi^{\prime}(0)=0$. Then the function
$\psi^{\prime}$ is cylindric if and only if $\psi$ is cylindric;
in this case we have
\[
(\mathbb{L}-1)\int_{\mathbb{P} X} \psi d \chi_g = \int_{X}
\psi^{\prime} d \chi_g
\]
(cf. \cite[(2.7)]{demo}).
}

We define now the generalised Poincaré series of the
projectivisation of the fractional ideal $\bgo$:

\textbf{(4.10)~Definition:} 
{\sl 
The \emph{generalised Poincar\'e series} of a multi-index
filtration given by the ideals $J(\underline{v})$ is the integral
\[
P_g (\bgo,\underline{t}, \mathbb{L}):=\int_{\mathbb{P} \bgo}
\underline{t}^{\underline{v}(z)} d \chi_g \in K_0
(\Nu_{k})_{(\mathbb{L})} [\![t_1, \ldots, t_r]\!],
\]
where $\underline{t}^{\underline{v}(z)}:=t_1^{v_1 (z)} \cdot
\ldots \cdot t_r^{v_r(z)}$ is considered as a (cylindric) function
on $\mathbb{P}\oo$ with values in $\mathds{Z}[\![t_1, \ldots,
t_r]\!]$ (the vector $\underline{v}(z)$ is supposed to be
$\underline{0}$ as soon as $v_i (z)=\infty$ for at least one $i
\in \{1, \ldots , r \}$).
}

\textbf{(4.11)~Remark:} 
{\sl 
Notice that if $\bgo=\oo$, then
$P_g(\oo,\underline{t},\mathbb{L})$ is the generalised Poincaré
series of a filtration $J(\underline{v})$ over the
projectivization of the ring $\oo$ introduced in \cite[Section 2,
p.~198]{cadegu11} for the case of the ring $\oo_{V,0}$ of
functions on a germ $(V,0)$ of a complex analytic variety.
}

\textbf{(4.12)} Let us define the \emph{degree} of a fractional $\oo$-ideal by the following two properties: (i) $\deg (\oo):=0$; (ii) for every two fractional $\oo$--ideals $\ago$, $\bgo$, one has $\deg (\ago) - \deg (\bgo) =  \dim (\ago/\bgo)$ whenever $\ago \supseteq \bgo$. Now, if we define
\[
L_g(\bgo, \underline{t}, \mathbb{L}):= \sum_{\underline{v} \in
\mathbb{Z}^r} \left ( \mathbb{L}^{\deg (J^{\bgo}(\underline{v}))}
- \mathbb{L}^{\deg (J^{\bgo}(\underline{v}+\underline{1}))} \right
) \cdot \underline{t}^{\underline{v}},
\]
then we get
\medskip

\textbf{(4.13)~Lemma:} 
{\sl 
\[
P_g (\bgo, \underline{t}, \mathbb{L}) =
\frac{\prod_{i=1}^{r}(t_i-1) L_g (\oo,\bgo; \underline{t})}{t_1
\cdot \ldots \cdot t_r-1}.
\]
}

\dem~ The result may be proved in much the same way as in
\cite[Proposition 2]{cadegu11}. \qed
\medskip

We describe now the functional equations for the series
$P_g(\bgo,\underline{t},\mathbb{L})$. First of all, we state the
following two results, due to St\"ohr (see \cite{stohr} and
\cite{stohr3}). We include the proofs by the sack of completeness.

\textbf{(4.14)~Lemma:} 
{\sl 
Let $\ago$, $\mathfrak{b}$ be fractional ideals of $\oo$ such that $\ago \supseteq \bgo$. We have
\[
\lambda (\bgo^{\ast} \cap \ago/\bgo \cap \ago^{\ast}) = \lambda (\ago / \bgo).
\]
}

\dem~  By definition of the ideal $\cgo $ we have $\cgo : \bgo^{\ast} = \bgo$, hence $\bgo^{\ast} \cap \ago= (\cgo : \bgo) \cap (\cgo : \ago^{\ast}) =\cgo :(\bgo + \ago^{\ast})$ and 
$\deg (\bgo^{\ast} \cap \ago)\overset{\star}{=} \deg( \cgo ) + \deg (\bgo + \ago^{\ast}) $. From the isomorphisms
\[
\bgo + \ago^{\ast}/\oo \cong \bgo + \ago^{\ast}/\ago^{\ast} + \ago^{\ast}/\oo \cong \bgo / \bgo \cap \ago^{\ast} +\ago^{\ast}/\oo 
\]
it follows that $\deg (\bgo + \ago^{\ast}) = \deg (\bgo) - \deg (\bgo \cap \ago^{\ast}) + \deg (\cgo) - \deg (\ago) $, and by $\star$ we get
\[
\deg (\bgo^{\ast} \cap \ago) = \deg (\bgo \cap \ago^{\ast}) + \deg (\ago) - \deg (\bgo),
\]
which proves the statement. \qed

\textbf{(4.15)~Lemma:} 
{\sl 
Let $\ago$ be a fractional ideal of $\ooo$. The following assertions hold:
\begin{itemize}
\item[(a)] Let $t_i$ be a generator of the ideal $\mm_i$ for
every $i \in \{1, \ldots ,r\}$. The fractional ideals $\ago$ of
$\oo$ are of the form $ \underline{t}^{-\underline{n}} \cdot \mathfrak{b}$,
where $\underline{t}^{\underline{n}}:=t_1^{n_1} \cdot \ldots \cdot
t_r^{n_r}$ for some $\underline{n}=(n_1, \ldots , n_r) \in
\mathds{Z}^r$ and being $\mathfrak{b}$ a fractional ideal of $\oo$
such that $\mathfrak{b} \cdot \overline{\oo} = \overline{\oo}$.

\item[(b)] There exists some $\underline{v} \in \mathds{Z}^r$ such
that $\ago = J^{\mathcal{K}}(\underline{v})$ and
$\ago^{\ast}=J^{\mathcal{K}}(-\underline{v})$.

\item[(c)] For some $\underline{v} \in \mathds{Z}^r$, we have
$\deg (\ago) = \delta - \underline{v} \cdot \underline{d}$.
\end{itemize}
}

\dem~\text{(a)} Since $\ago$ is a fractional ideal of $\oo$, $\ago
\overline{\oo}$ is a fractional ideal of $\overline{\oo}$. By
Remark \ref{rem:cero}, we have that $\ago \overline{\oo}$ must be
of the form $\mm_1^{n_1} \cdot \ldots \cdot \mm_r^{n_r}$ for some
$\underline{n} \in \mathds{Z}^r$. That is,
\[
\ago \overline{\oo} =  t_1^{n_1} \overline{\oo} \cdot \ldots \cdot t_r^{n_r} \overline{\oo}  = \underline{t}^{\underline{n}} \cdot \overline{\oo}, 
\]
i.e., $\underline{t}^{-\underline{n}} \cdot \ago \cdot
\overline{\oo} = \overline{\oo}$ for some $\underline{n} \in
\mathds{Z}^{r}$. Then, it suffices to take
$\mathfrak{b}=\underline{t}^{-\underline{n}}\cdot \ago$ and the claim
follows.

\text{(b)} Let $\cgo$ be a canonical ideal of $\oo$. We know that
\begin{eqnarray}
\deg (\cgo)& = & \deg (\cgo : \oo)  \nonumber \\
& = & \deg (\cgo : \ooo) + \dim \left ( \ooo / \oo \right )  \nonumber \\
& = & \deg (\cgo : \ooo) + \delta. \nonumber
\end{eqnarray}
If we multiply the ideal $\cgo$ by a convenient element of
$\mathcal{K}$, then we may assume $\cgo : \ooo = \ooo$. The ideal
$\cgo : \ooo$ is fractional, hence by (a) there is some
$\underline{v} \in \mathds{Z}^r$ such that $\cgo : \ooo =
\underline{t}^{\underline{v}} \cdot \ooo$, which is equivalent to $\ooo =
(\underline{t}^{-\underline{v}} \cdot \cgo) : \ooo$. Then we have
\[
\ago^{\ast}\cgo : \ago  = \cgo : J^{\mathcal{K}}(\underline{v}) = \cgo : (\underline{t}^{\underline{v}} \cdot \ooo) = \underline{t}^{-\underline{v}}(\cgo : \ooo) = \underline{t}^{-\underline{v}} \cdot \ooo  = J^{\mathcal{K}}(-\underline{v}). 
\]


\text{(c)} By $(a)$, there is some $\underline{v} \in
\mathds{Z}^r$ such that $\ago=J^{\mathcal{K}}(\underline{v})$.
Since $d_i = \dim \left ( \ooo / \mm_i \right )$ for all $1 \le i
\le r$, by the Chinese Remainder Theorem we have
\begin{eqnarray}
\dim \left ( \ooo / J^{\mathcal{K}}(\underline{v}) \right ) & = & \deg (\ooo) - \deg (J^{\mathcal{K}}(\underline{v})) \nonumber \\
 & = & \sum_{i=1}^{r} d_i \cdot v_i \nonumber \\
 & = & \underline{d} \cdot \underline{v}. \nonumber
\end{eqnarray}
Since $\delta = \dim \left ( \ooo / \oo \right ) = \deg (\ooo)$,
we have
\begin{eqnarray}
\deg (J^{\mathcal{K}}(\underline{v}))& = & \deg (\ooo) - \underline{d} \cdot \underline{v} \nonumber \\
& = & \delta - \underline{d} \cdot \underline{v}.\nonumber
\end{eqnarray}
 \qed

Next proposition relates the degree of the ideal
$J^{\bgo^{\ast}}(\underline{v})$ and the value $\gamma^{\bgo}$.

\textbf{(4.16)~Proposition:} 
{\sl 
For every $\underline{v} \in \mathds{Z}^r$, we have
\begin{eqnarray}
\deg \left ( J^{\bgo^{\ast}}(\underline{v}) \right ) & = & \deg \left ( J^{\bgo}(\gamma^{\bgo} -\underline{v})  \right ) + \dim \left ( \bgo / (\bgo : \ooo) \right ) - \underline{v} \cdot \underline{d}  \nonumber \\
 & = & \deg \left ( J^{\bgo}(\gamma^{\bgo} -\underline{v})  \right ) + \dim \left ( \bgo^{\ast} \cdot \ooo / \bgo^{\ast} \right )  - \underline{v} \cdot \underline{d}  \nonumber
\end{eqnarray}

where $\underline{v} \cdot \underline{d}:=v_1 d_1 + \ldots + v_r
d_r$.
}

\dem~ The second equality holds by Lemma \text{(2.11)}.
Moreover, since $J^{\bgo}(\underline{v})=\bgo \cap
J^{\mathcal{K}}(\underline{v})$ for $\underline{v} \in
\mathds{Z}^r$, by Lemma \text{(4.14)} and Lemma
\text{(4.15)} it follows that
\[
\deg \left ( \bgo^{\ast} \cap J^{\mathcal{K}}(\underline{v})
\right ) = \dim \left ( \bgo / \bgo : \ooo \right ) + \deg \left (
J^{\mathcal{K}}(-\underline{v}) \cdot \overline{\oo^{\ast}} \cap
\bgo \right ) - \underline{v} \cdot \underline{d}.
\]
The definition of $\gamma^{\bgo}$ allows us now to conclude. \qed

As a consequence we obtained the following result due to Moyano and Z\'u\~niga (\cite[Lemma 9]{mozu})

\textbf{(4.17)~Corollary:} 
{\sl 
The ring $\oo$ is Gorenstein if and only if
\[
\ell (\gamma - \underline{v}) - \ell (\underline{v}) = \delta -
\underline{v} \cdot \underline{d}
\]
for every $\underline{v} \in \mathds{Z}^r$.
}

\dem~ It is just to apply Proposition \text{(4.16)} to
$\bgo=\oo$. Notice that we have $\oo=\oo^{\ast}$ (cf.
\text{(2.7)}) because $\oo$ is Gorenstein. \qed
\medskip

\textbf{(4.18)~Remark:} 
{\sl 
Notice that $\delta - d = \ell(\gamma - \underline{1}) - \rho$, if
the ring is Gorenstein. It follows from Corollary
\text{(4.17)}, because $\delta - \underline{d} \cdot
\underline{1} = \ell (\gamma - \underline{1}) -
\ell(\underline{1})$ and $\ell (\underline{1})=\rho$.
}

Proposition \text{(4.16)} allows us to describe the functional equations for the generalised Poincar\'e series:

\textbf{(4.19)~Theorem:} 
{\sl 
\[
L_g (\bgo, \mathbb{L}^{d_1}t_1, \ldots , \mathbb{L}^{d_r}t_r,
\mathbb{L})= \mathbb{L}^{\dim \left ( \bgo^{\ast} \cdot \ooo /
\bgo^{\ast} \right )-d} \cdot
\underline{t}^{\gamma^{\bgo}-\underline{1}} \cdot L_g (
\bgo^{\ast};\underline{t}^{-1}, \mathbb{L}).
\]
}

\dem~ Let $A (\bgo,\underline{t}):= \sum_{\underline{v} \in
\mathds{Z}^r} \mathbb{L}^{\deg \left ( J^{\bgo}(\underline{v})
\right )} \cdot \underline{t}^{\underline{v}}$. Then
\begin{eqnarray}
 L_g(\bgo, \underline{t}, \mathbb{L}) & = & \sum_{\underline{v} \in \mathds{Z}^r} \left ( \mathbb{L}^{\deg \left ( J^{\bgo}(\underline{v})
\right )} - \mathbb{L}^{\deg \left ( J^{\bgo}(\underline{v} +
\underline{1})
\right )} \right ) \cdot \underline{t}^{\underline{v}} \nonumber \\
  & = & (1-\underline{t}^{-1}) \cdot A(\bgo, \underline{t}). \nonumber \label{eqn:uno}
\end{eqnarray}

We apply now Proposition \text{(4.16)} to obtain

\begin{equation} \label{eqn:dos}
A(\bgo,\mathbb{L}^{d_1}t_1, \ldots , \mathbb{L}^{d_r}t_r) =
\mathbb{L}^{\dim \left ( \bgo^{\ast} \cdot \ooo / \bgo^{\ast}
\right )} \cdot \underline{t}^{\gamma^{\bgo}} \cdot A
(\bgo^{\ast},t_1^{-1}, \ldots , t_r^{-1}).\nonumber
\end{equation}

Moreover, taking the inverse of $\underline{t}$, we have
\[
L_g (\bgo^{\ast},\underline{t}^{-1}, \mathbb{L})  =
(1-\underline{t}) \cdot A(\bgo^{\ast},\underline{t}^{-1}).
\nonumber \tag{$\dag$}
\]

Therefore
\begin{align*}
L_g (\oo,\bgo,\mathbb{L}^{d_1}t_1, \ldots , \mathbb{L}^{d_r}t_r)  = & \frac{\mathbb{L}^d \cdot \underline{t}-1}{\mathbb{L}^d \cdot \underline{t}} \cdot A (\bgo,\mathbb{L}^{d_1}t_1, \ldots , \mathbb{L}^{d_r}t_r) \\
  = & \frac{\mathbb{L}^d \cdot \underline{t}-1}{\mathbb{L}^d \cdot \underline{t}} \cdot \mathbb{L}^{\dim \left ( \bgo^{\ast} \cdot \ooo / \bgo^{\ast} \right )} \cdot \underline{t}^{\gamma^{\bgo}} \cdot A(\bgo^{\ast},t_1^{-1}, \ldots , t_r^{-1}) \\
  \overset{=}{(\dag)} & \frac{\mathbb{L}^d \cdot \underline{t}-1}{\mathbb{L}^d \cdot \underline{t}} \cdot \frac{1}{1-t} \cdot \mathbb{L}^{\dim \left ( \bgo^{\ast} \cdot \ooo / \bgo^{\ast} \right )} \cdot \underline{t}^{\gamma^{\bgo}} \cdot L_g (\bgo^{\ast},\underline{t}^{-1},\mathbb{L}) \\
  = & \frac{\mathbb{L}^d \cdot \underline{t}-1}{1-\underline{t}} \cdot \mathbb{L}^{\dim \left ( \bgo^{\ast} \cdot \ooo / \bgo^{\ast} \right ) - d} \cdot \underline{t}^{\gamma^{\bgo} - \underline{1}} \cdot L_g(\bgo^{\ast},\underline{t}^{-1}, \mathbb{L}).
 \end{align*}
\qed

\textbf{(4.20)~Corollary:} 
{\sl 
\[
P_g (\bgo,\mathbb{L}^{d_1}t_1, \ldots , \mathbb{L}^{d_r}t_r,
\mathbb{L}) = \mathbb{L}^{\dim \left ( \bgo^{\ast} \cdot \ooo /
\bgo^{\ast} \right ) - d} \cdot
\underline{t}^{\gamma^{\bgo}-\underline{1}} \cdot
\frac{\prod_{i=1}^{r}(1-\mathbb{L}^{d_i}t_i)}{\prod_{i=1}^{r}(t_i-1)}
\cdot P_g (\bgo^{\ast},\underline{t}^{-1},\mathbb{L}).
\]
}

\textbf{(4.21)~Corollary:} 
{\sl 
If $\bgo=\oo$, then we have
\[
P_g (\oo,\mathbb{L}^{d_1}t_1, \ldots ,
\mathbb{L}^{d_r}t_r,\mathbb{L}) = \mathbb{L}^{\dim \left (
\oo^{\ast} \cdot \ooo / \oo^{\ast} \right ) - d} \cdot
\underline{t}^{\gamma-\underline{1}} \cdot
\frac{\prod_{i=1}^{r}(1-\mathbb{L}^{d_i}t_i)}{\prod_{i=1}^{r}(t_i-1)}
\cdot P_g (\oo^{\ast},\underline{t}^{-1},\mathbb{L}).
\]
Furthermore, if $\oo$ is Gorenstein, then we obtain
\[
P_g (\oo,\mathbb{L}^{d_1}t_1, \ldots ,
\mathbb{L}^{d_r}t_r,\mathbb{L}) = \mathbb{L}^{\delta - d} \cdot
\underline{t}^{\gamma-\underline{1}} \cdot
\frac{\prod_{i=1}^{r}(1-\mathbb{L}^{d_i}t_i)}{\prod_{i=1}^{r}(t_i-1)}
\cdot P_g (\oo,\underline{t}^{-1},\mathbb{L}).
\]
}

\dem~ It is a straight consequence of Corollary \text{(4.20)} and
Corollary \text{(4.17)}. \qed

\section{\label{functional2}Extended generalised value ideal Poincaré series}

Let $\oo$ be a one-dimensional Cohen-Macaulay local Noetherian
ring having a perfect coefficient field $K$. Campillo, Delgado and
Gusein-Zade introduced in \cite{extended} the notion of extended
semigroup of a germ of complex plane curve singularity. We want
now to define the concept of extended value ideal of a fractional
ideal $\bgo$. Let us preserve notations as in \text{(2.1)}. Remember that the
ideal $\mm_j V_j$ is regular maximal of $V_j$ and
$V_j=\ooo_{\mm_j}$ for every $j \in \{1, \ldots , r \}$.
\medskip

If $K_j$ is a coefficient field of $V_j$ and $t_j$ is an
indeterminate over $K_j$, then one can identify $V_j \cong K_j
[\![t_j]\!]$ and $v_j$ with the order function respect to $t_j$ in
$K_j [\![t_j]\!]$ for every $j \in \{1, \ldots , r \}$. Thus
\[
\oo \subset K_1 [\![t_1]\!] \cap \ldots \cap K_r [\![t_r]\!] =
\ooo.
\]

Since $V_j$ is an $\oo$-module of finite type, the field
extensions $\oo / \mm \hookrightarrow \ooo / \mm_j$ are finite for
every $j \in \{1, \ldots r \}$. Furthermore, as $\oo / \mm$ is
assumed to be perfect, every such a extension is separable and
therefore, for every coefficient field $K$ of $\oo$ there exists a
unique coefficient field $K_j$ of $V_j$ with $K \subset K_j$ which
is isomorphic to $\ooo / \mm_j$ for every $j \in \{1, \ldots , r
\}$.
\medskip

Let us consider the vector spaces $C^{\bgo}(\underline{v},i)=
J^{\bgo}(\underline{v})/J^{\bgo}(\underline{v}+\underline{1}_{\{ i
\}})$ for every $i \in \{1, \ldots , r \}$ and the map:
\begin{displaymath}
\begin{array}{lccc}
j_{\underline{v}}: & J^{\bgo}(\underline{v}) & \longrightarrow & C^{\bgo}(\underline{v},1) \times \ldots \times C^{\bgo}(\underline{v},r)  \\
& z & \mapsto & \left ( j_1 (z), \ldots , j_r (z) \right
)=:j_{\underline{v}}(z) .
\end{array}
\end{displaymath}
We can identify
$\mathrm{Im~}j_{\underline{v}} \cong C^{\bgo}(\underline{v})=
J^{\bgo}(\underline{v})/J^{\bgo}(\underline{v}+\underline{1})$ and
define the set
\[
F^{\bgo}_{\underline{v}}:= C^{\bgo}(\underline{v}) \cap \left (
(C^{\bgo}(\underline{v},1) \setminus \{0\}) \times \ldots \times
(C^{\bgo}(\underline{v},r) \setminus \{0\}) \right ).
\]
\textbf{(5.1)~Lemma:} 
{\sl 
\[
F^{\bgo}_{\underline{v}} = C^{\bgo}(\underline{v}) \cap \left (
K_1^{\ast} \times \ldots \times K_r^{\ast} \right ).
\]
}

\dem~It is enough to define an isomorphism $\varphi_{\underline{v}}: C^{\bgo}(\underline{v},1)
\setminus \{0 \} \to K_1^{\ast}$. Let $z \in
\bgo \setminus \{0 \}$ with $v_1(z)=v_1$. We have that $z = a_1(z)
t_1^{v_1(z)}$ with $a_1(z) \in K_1^{\ast}$, thus $\varphi_{\underline{v}}$ can be
defined by $z \mapsto a_1(z)$. \qed

\medskip

\textbf{(5.2)} For every $\underline{v} \in \mathds{Z}^r$ we write
\[
F^{\bgo}_{\underline{v}}= \left ( J^{\bgo}(\underline{v}) /
J^{\bgo}(\underline{v}+\underline{1}) \right ) \setminus
\bigcup_{i=1}^{r} \left ( J^{\bgo}(\underline{v} +
\underline{1}_{\{i \}}) / J^{\bgo}(\underline{v}+\underline{1})
\right ),
\]
i.e., $F^{\bgo}_{\underline{v}}$ is the complement to an
arrangement of vector subspaces in a vector space (it is not a
vector subspace itself). Notice that this is precisely the space
(2) in \text{(4.2)}.

\textbf{(5.3)~Definition:} 
{\sl 
For every fractional ideal $\bgo$ of $\oo$ we define the set
$\widehat{S}(\bgo)$ to be the union of the subspaces
$F^{\bgo}_{\underline{v}}$ for all $\underline{v} \in
\mathds{Z}^r$. The spaces $F^{\bgo}_{\underline{v}}$ are called
\emph{fibres} of $\widehat{S}(\bgo)$.
}

Notice that, if $\bgo=\oo$, then $\widehat{S}(\oo)$ is the
extended semigroup of the ring $\oo$ (see \cite{extended} for
further details).

\textbf{(5.4)} The group $K^{\ast}$ of non-zero elements of $K$ acts freely
on $\mathds{Z}^r \times (K_1^{\ast} \times \ldots \times
K_r^{\ast})$ (by multiplication of all coordinates in $K_1^{\ast}
\times \ldots \times K_r^{\ast}$). The corresponding factor space
$\mathds{Z}^{r} \times (K_1^{\ast} \times \ldots \times
K_r^{\ast}) / K^{\ast} = \mathds{Z}^r \times \pro (K_1^{\ast}
\times \ldots \times K_r^{\ast})= \sum_{\underline{v} \in
\mathds{Z}^r} \pro (K_1^{\ast} \times \ldots \times K_r^{\ast})
\underline{t}^{\underline{v}}$ has the natural structure of
semigroup. The set $\widehat{S}(\bgo) \subset \mathds{Z}^r \times
(K_1^{\ast} \times \ldots \times K_r^{\ast})$ is invariant with
respect to the $K^{\ast}$-action. The factor space
\[
\pro \widehat{S}(\bgo) = \widehat{S}(\bgo) / K^{\ast}
\]
is called the projectivisation of $\widehat{S}(\bgo)$ (it is also
a graded $S(\oo)$-module in a natural sense).
\medskip

From the previous definitions, the projectivisation of
$\widehat{S}(\bgo)$ can be described as
\[
\pro \widehat{S}(\bgo) = \sum_{\underline{v} \in \mathds{Z}^r}
\pro F^{\bgo}_{\underline{v}} \cdot \underline{t}^{\underline{v}},
\]
where $\pro F^{\bgo}_{\underline{v}}= F^{\bgo}_{\underline{v}} /
K^{\ast}$ is the projectivisation of the fibre
$F^{\bgo}_{\underline{v}}$. For $\underline{v} \in
\widehat{S}(\bgo)$, the space $\pro F^{\bgo}_{\underline{v}}$ is
the complement to an arrangement of projective hyperplanes in a
$\left (\dim_{K} \left ( J^{\bgo}(\underline{v}) /
J^{\bgo}(\underline{v}+\underline{1}) \right ) - 1
\right)$-dimensional projective space $\pro \left (
J^{\bgo}(\underline{v})/J^{\bgo}(\underline{v}+\underline{1})
\right )$.
\medskip

Set the Laurent series
\[
\chi_g (\pro \widehat{S}(\bgo)):= \sum_{\underline{v} \in
\mathds{Z}^r} \chi_g \left ( \pro F^{\bgo}_{\underline{v}} \right
) \cdot \underline{t}^{\underline{v}}.
\]

On the other hand, we also have the extended generalised value ideal Poincaré
series of a filtration $\{ J^{\bgo}(\underline{v}) \}$ defined by
$\underline{v}(z)= (v_1 (z), \ldots , v_r(z))$, for $z \in \bgo$:
\[
\widehat{P}_g (\bgo, \underline{t}, \mathbb{L}, \{v_i \}) :=
\int_{\pro \widehat{S}(\bgo)} \underline{t}^{\underline{v}(z)} d
\chi_g
\]
(we will use $\widehat{P}_g (\bgo, \underline{t}, \mathbb{L})$
instead of $\widehat{P}_g (\bgo, \underline{t}, \mathbb{L}, \{v_i
\})$ when the filtration is clear from the context). Notice that if
$\mathfrak{b}=\oo$, then $\widehat{P}_g (\oo, \underline{t},
\mathbb{L})$ coincides with the generalised semigroup Poincaré
series defined in \cite[p. 507]{cadegu11}.
\medskip

All projectivisations $\pro F^{\bgo}_{\underline{v}}$ of the
fibres $F^{\bgo}_{\underline{v}}$ (i.e., all connected components
of $\pro \widehat{S}(\bgo)$) are complements to arrangements of
projective subspaces in finite dimensional projective spaces. We
define
\[
\widehat{L}_g(\bgo, \underline{t}, \mathbb{L}):=
\sum_{\underline{v} \in \mathds{Z}^r} \left [ \pro
(J^{\bgo}(\underline{v}) / J^{\bgo}(\underline{v}+\underline{1}))
\right ] \cdot \underline{t}^{\underline{v}}.
\]

\textbf{(5.5)~Proposition:} 
{\sl 
\[
\chi_g (\pro \widehat{S}(\bgo)) = \widehat{P}_g
(\bgo,\underline{t},\mathbb{L}) = \frac{\widehat{L}_g(\bgo,
\underline{t}, \mathbb{L}) \cdot \prod_{i=1}^{r} (t_i - 1)}{t_1
\cdot \ldots \cdot t_r - 1}.
\]
}

\dem~Let $\underline{w} \in \mathds{Z}^r$ and set
$L_{I}:= \{ (a_1, \ldots , a_r) \in K^r \mid a_i=0 \mathrm{~for~}
i \in I \}$. Then
\begin{align*}
\chi_g (\pro F^{\bgo}_{\underline{v}}) = & \chi_g \left ( \pro J^{\bgo}(\underline{v})/J^{\bgo}(\underline{v}+\underline{1}) \right ) - \chi_g \left ( \bigcup_{i=1}^{r} \pro \left ( J^{\bgo}(\underline{v})/J^{\bgo}(\underline{v}+\underline{1}) \cap L_{\{ i\}} \right ) \right )     \\
    = & \chi_g \left (\pro J^{\bgo}(\underline{v})/J^{\bgo}(\underline{v}+\underline{1}) \right ) - \sum_{\substack{I \subset I_0 \\ I \ne \varnothing}} (-1)^{\sharp I - 1} \chi_g \left ( \pro \left ( J^{\bgo}(\underline{v})/J^{\bgo}(\underline{v}+\underline{1}) \cap L_I \right ) \right )  \\
    = & \sum_{I \subset I_0} (-1)^{\sharp I} \chi_g \left (\pro \left( J^{\bgo}(\underline{v})/J^{\bgo}(\underline{v}+\underline{1}) \cap L_I \right )  \right )    \\
    = & \sum_{I \subset I_0} (-1)^{\sharp I} \left [ \pro \left ( J^{\bgo}(\underline{v})/J^{\bgo}(\underline{v}+\underline{1}) \cap L_I \right ) \right ]   \\
    = & \sum_{I \subset I_0} (-1)^{\sharp I} \left [ \pro \left ( J^{\bgo}(\underline{v}+\underline{1}_{I})/J^{\bgo}(\underline{v}+\underline{1}) \right )\right ].
\end{align*}

Therefore
\begin{align*}
(t_1 \cdot \ldots \cdot t_r -1)\chi_g (\pro F^{\bgo}_{\underline{v}})  = & \sum_{\underline{v} \in \mathds{Z}^r} \sum_{I \subset I_0} (-1)^{\sharp I} \left [ \pro \left ( J^{\bgo}(\underline{v}+\underline{1}_{I}-\underline{1})/J^{\bgo}(\underline{v}) \right )\right ] \cdot \underline{t}^{\underline{v}}- \\
    - & \sum_{\underline{v} \in \mathds{Z}^r} \sum_{I \subset I_0} (-1)^{\sharp I} \left [ \pro \left ( J^{\bgo}(\underline{v}+\underline{1}_{I})/J^{\bgo}(\underline{v}+\underline{1}) \right )\right ] \cdot
    \underline{t}^{\underline{v}}  \\
    = & \sum_{\underline{v} \in \mathds{Z}^r} \sum_{I \subset I_0} (-1)^{\sharp I} \left [ \pro \left ( J^{\bgo}(\underline{v}+\underline{1}_{I}-\underline{1})/J^{\bgo}(\underline{v}+\underline{1}_{I}) \right )\right ] \cdot \underline{t}^{\underline{v}}.
    \ \ \ \ (\ast)
\end{align*}
The coefficient of $\underline{t}^{\underline{v}}$ in the
polynomial
\[
\left ( \sum_{\underline{v} \in \mathds{Z}^r} \left [ \pro \left (
J^{\bgo}(\underline{v})/J^{\bgo}(\underline{v}+\underline{1})
\right ) \right ] \cdot \underline{t}^{\underline{v}} \right )
\cdot \prod_{i=1}^{r} (t_i -1)
\]
is equal to
\[
\sum_{I \subset I_0} (-1)^{\sharp I} \left [ \pro \left (
J^{\bgo}(\underline{v}-\underline{1}+\underline{1}_{I})/J^{\bgo}(\underline{v}+\underline{1}_{I})
\right ) \right ],
\]
and the latter formula coincides with $(\ast)$. \qed
\medskip

\textbf{(5.6)~Remark:} 
{\sl 
Since $\mathbb{L}-1$ is invertible in $K_0 (\Nu_{k})_{(\mathbb{L})}$, the extended generalised value ideal Poincar\'e series can be
rewritten as
\[
\widehat{P}_g (\bgo,\underline{t}, \mathbb{L}) =
\frac{\prod_{i=1}^{r}(t_i - 1)}{t_1 \cdot \ldots \cdot t_r-1}
\cdot \sum_{\underline{v} \in \mathds{Z}^r}
\frac{\mathbb{L}^{c^{\bgo}(\underline{v})}-1}{\mathbb{L}-1} \cdot
\underline{t}^{\underline{v}}.
\]
}
\medskip

We generalise this functional equation by using the following
result:

\textbf{(5.7)~Proposition:} 
{\sl 
For every $\underline{v} \in \mathds{Z}^r$, we have
\[
c^{\bgo^{\ast}} (\underline{v}) = d - c^{\bgo}
(\gamma^{\bgo}-\underline{v}-\underline{1}).
\]
}

\dem~By means of Proposition \text{(4.16)} we deduce the equality:
\begin{eqnarray}
c^{\bgo^{\ast}}(\underline{v}) & = & \deg \left (J^{\bgo^{\ast}} (\underline{v}) \right ) - \deg \left ( J^{\bgo^{\ast}} (\underline{v}+\underline{1}) \right )\nonumber \\
 & = & \dim \left ( \bgo / \bgo:\ooo \right ) - \underline{v} \cdot \underline{d}+\deg
 \left ( J^{\bgo} (\gamma^{\bgo}-\underline{v}) \right ) \nonumber \\
 & & - \dim \left ( \bgo / \bgo:\ooo \right ) - (\underline{v}+ \underline{1}) \cdot \underline{d}-\deg
 \left ( J^{\bgo} (\gamma^{\bgo}-\underline{v}-\underline{1}) \right )
 \nonumber \\
& = & d - c^{\bgo} (\gamma^{\bgo}-\underline{v}-\underline{1}).
\nonumber
\end{eqnarray}
\qed

We want to describe functional equations for the series
$\widehat{L}_g(\bgo, \underline{t}, \mathbb{L})$ and the Poincaré
series $\widehat{P}_g(\bgo,\underline{t}, \mathbb{L})$.

\textbf{(5.8)~Theorem:} 
{\sl 
\[
\underline{t}^{\gamma^{\bgo}-\underline{1}} \widehat{L}_g (\bgo,
\underline{t}^{-1}, \mathbb{L}) = - \mathbb{L}^{d-1} \widehat{L}_g
(\bgo^{\ast},\underline{t},\mathbb{L}^{-1}).
\]
\[
\underline{t}^{\gamma^{\bgo}-\underline{1}} \widehat{P}_g (\bgo,
\underline{t}^{-1}, \mathbb{L}) = (-1)^r \mathbb{L}^{d-1}
\widehat{P}_g (\bgo^{\ast},\underline{t},\mathbb{L}^{-1}).
\]
}

\dem~By Remark \text{(5.6)} we have
$\widehat{L}_g(\bgo,\underline{t},\mathbb{L})=\sum_{\underline{v}
\in \mathds{Z}^r}
\frac{\mathbb{L}^{c(\underline{v})}-1}{\mathbb{L}-1} \cdot
\underline{t}^{\underline{v}}$. It suffices to take
$\widetilde{L}_g(\bgo,\underline{t},\mathbb{L})=\sum_{\underline{v}
\in \mathds{Z}^r}
\frac{\mathbb{L}^{c(\underline{v})}}{\mathbb{L}-1} \cdot
\underline{t}^{\underline{v}}$. By Proposition \text{(5.7)} it
holds
$c^{\bgo^{\ast}}(\underline{v})+c^{\bgo}(\gamma-\underline{v}-\underline{1})=d$
and we have
\[
\underline{t}^{\gamma^{\ast}-\underline{1}} \cdot
\widetilde{L}_g(\bgo,\underline{t}^{-1},\mathbb{L})=\mathbb{L}^d
\cdot \sum_{\underline{v} \in \mathds{Z}^r}
\frac{\mathbb{L}^{-c^{\bgo^{\ast}}(\underline{v})}}{\mathbb{L}-1}
\cdot \underline{t}^{\underline{v}}.
\]
On the other hand, we have
\begin{eqnarray}
\widetilde{L}_g(\bgo^{\ast},\underline{t}, \mathbb{L}^{-1})& = & \sum_{\underline{v} \in \mathds{Z}^r} \frac{\mathbb{L}^{-c^{\bgo^{\ast}}(\underline{v})}}{\mathbb{L}^{-1}-1} \cdot \underline{t}^{\underline{v}}  \nonumber \\
  & = & -\mathbb{L} \cdot \sum_{\underline{v} \in \mathds{Z}^r}
  \frac{\mathbb{L}^{-c^{\bgo^{\ast}}(\underline{v})}}{\mathbb{L}-1} \cdot \underline{t}^{\underline{v}}.
  \nonumber
\end{eqnarray}
Hence
\[
-\mathbb{L}^{d-1} \cdot
\widetilde{L}_g(\bgo^{\ast},\underline{t},\mathbb{L}^{-1})=\underline{t}^{\gamma^{\bgo}-\underline{1}}
\cdot \widetilde{L}_g(\bgo,\underline{t}^{-1},\mathbb{L}).
\]
Taking now into account that
\[
\widehat{P}_g(\bgo,\underline{t},\mathbb{L})=\frac{(t_1-1) \cdot
\ldots \cdot (t_r-1)}{t_1 \cdot \ldots \cdot t_r - 1} \cdot
\widetilde{L}_g(\bgo,\underline{t},\mathbb{L}),
\]
and the fact that
\[
(-1)^{r-1} \cdot \frac{(t_1-1) \cdot \ldots \cdot (t_r-1)}{t_1
\cdot \ldots \cdot t_r -1} = \frac{(1-t_1) \cdot \ldots \cdot
(1-t_r)}{1-t_1 \cdot \ldots \cdot t_r},
\]
we get
\begin{eqnarray}
\widehat{P}_g(\bgo,\underline{t}^{-1},\mathbb{L}) & = & \frac{(t_1^{-1}-1) \cdot \ldots \cdot (t_r^{-1}-1)}{t_1^{-1}\cdot \ldots \cdot t_r^{-1}-1} \cdot \widetilde{L}_g(\bgo,\underline{t}^{-1},\mathbb{L})     \nonumber \\
  & = & \frac{(1-t_1) \cdot \ldots \cdot (1-t_r)}{1-t_1 \cdot \ldots \cdot
  t_r} \cdot \widetilde{L}_g(\bgo,\underline{t}^{-1},\mathbb{L})
  \nonumber\\
  & = & (-1)^{r-1} \cdot \frac{(t_1-1) \cdot \ldots \cdot (t_r-1)}{t_1 \cdot \ldots \cdot t_r - 1} \cdot \widetilde{L}_g(\bgo,\underline{t}^{-1},\mathbb{L}).\nonumber
\end{eqnarray}
Since
\begin{eqnarray}
\widehat{P}_g(\bgo^{\ast},\underline{t},\mathbb{L}^{-1}) & = &
\frac{(t_1-1) \cdot \ldots \cdot (t_r-1)}{t_1\cdot \ldots \cdot
t_r-1} \cdot
\widetilde{L}_g(\bgo^{\ast},\underline{t},\mathbb{L}^{-1}),
\nonumber
\end{eqnarray}
we have
\begin{eqnarray}
\underline{t}^{\gamma^{\ast}-\underline{1}} \cdot \widehat{P}_g(\bgo,\underline{t}^{-1},\mathbb{L}) & = & (-1)^{r-1} \cdot \frac{(t_1-1) \cdot \ldots \cdot (t_r-1)}{t_1 \cdot \ldots \cdot t_r - 1} \cdot \underline{t}^{\gamma^{\ast}-\underline{1}} \cdot \widetilde{L}_g(\bgo,\underline{t}^{-1},\mathbb{L}) \nonumber \\
  & = & (-1)^{r-1} \cdot \frac{(t_1-1) \cdot \ldots \cdot (t_r-1)}{t_1 \cdot \ldots \cdot t_r - 1} \cdot (-\mathbb{L})^{d-1} \cdot \widetilde{L}_g(\bgo^{\ast},\underline{t},\mathbb{L}^{-1}) \nonumber\\
  & = & (-1)^{r} \cdot \frac{(t_1-1) \cdot \ldots \cdot (t_r-1)}{t_1 \cdot \ldots \cdot t_r -
  1}\cdot
  \mathbb{L}^{d-1} \cdot
  \widetilde{L}_g(\bgo^{\ast},\underline{t},\mathbb{L}^{-1})\nonumber \\
  & = & (-1)^r \cdot \mathbb{L}^{d-1} \cdot \widehat{P}_g(\bgo^{\ast},\underline{t}, \mathbb{L}^{-1}),  \nonumber
\end{eqnarray}
and we are done. \qed

\end{document}